\documentclass[11pt]{amsart}
 \usepackage{amssymb,amsthm}

\usepackage[]{graphicx}
\usepackage{psfrag}

\everymath{\displaystyle}

\newtheorem{theo}{Theorem}[section]
\newtheorem{prop}[theo]{Proposition}

\newtheorem{lemma}[theo]{Lemma}
\newtheorem{teo}[theo]{Theorem}

\newtheorem{cor}[theo]{Corollary}
\theoremstyle{definition}
\newtheorem{defi}[theo]{Definition}

\newtheorem{conj}[theo]{Conjecture}
\newtheorem{remark}[theo]{Remark}

\newtheorem{question}[theo]{Question}

\newcommand{\fin}{\hfill$\square$}

\newcommand{\wt}{\widetilde}

\newcommand{\Cl}{\mbox{Cl}}

\newcommand{\sG}{\mathfrak G}
\newcommand{\sS}{\mathfrak S}
\newcommand{\sH}{\mathfrak H}

\newcommand{\sA}{\mathfrak A}
\newcommand{\sK}{\mathfrak K}

\newcommand{\sk}{\mathfrak k}

\newcommand{\ssl}{\mathfrak l}

\newcommand{\cN}{\mathcal N}

\newcommand{\DD}{\mathbb{D}}

\newcommand{\HH}{\mathbb{H}}

\newcommand{\NN}{\mathbb{N}}

\newcommand{\PP}{\mathbb{P}}

\newcommand{\RR}{\mathbb{R}}
\renewcommand{\SS}{\mathbb{S}}

\newcommand{\op}{\operatorname}

\DeclareMathOperator{\SO}{SO}
\DeclareMathOperator{\T}{T}
\DeclareMathOperator{\AdS}{AdS}
\DeclareMathOperator{\Ein}{Ein}

\newcommand{\uAdS}{\widetilde{\AdS}}

\newcommand{\ie}{\emph{ie. }}

\newcommand{\sca}[2]{\langle #1 | #2\rangle}         
 
\newcommand{\mcal}[1]{\ensuremath{\mathcal{#1}}}
\newcommand{\orth}{\bot}
\newcommand{\cY}{\mcal{Y}}

\newcommand{\cD}{\mcal{D}}
\newcommand{\cE}{\mcal{E}}
\newcommand{\cF}{\mcal{F}}

\newcommand{\cC}{\mcal{C}}

\begin{document}

\title[Quasi-Fuchsian AdS representations are Anosov]{Quasi-Fuchsian AdS representations are Anosov}



\author[T. Barbot]{Thierry Barbot}

\address{${}^\dagger$CNRS, UMPA, \'Ecole Normale Sup\'erieure de Lyon}
\email{Thierry.Barbot@umpa.ens-lyon.fr}


\begin{abstract}
Let $\Gamma$ be a cocompact lattice in $\SO(1,n)$. In \cite{part1} Q. M\'erigot
proved that if a representation $\rho: \Gamma \to \SO(2,n)$ is Anosov in the sense of
Labourie (\cite{labourie}), then it is \textit{quasi-Fuchsian}, \ie it is
faithfull, discrete, and preserves an acausal subset in the boundary of anti-de Sitter space.
In the present paper, we prove the reverse implication. It also includes:

-- A construction  of Dirichlet domains in the context of anti-de Sitter geometry,

-- A proof that spatially compact globally hyperbolic anti-de Sitter spacetimes with acausal
limit set admit locally $\op{CAT}(-1)$ Cauchy hypersurfaces.
\end{abstract}

\maketitle


\section{Introduction}
This paper is a complement to the preceding one by Q. M\'erigot \cite{part1}.
We will use all the notations, definitions and results therein. 
Let's just remind that a representation $\rho: \Gamma \to \op{SO}_{0}(2,n)$ is quasi-Fuchsian, or GHC-regular
if it is faithfull, dicrete and preserves an achronal topological $(n-1)$-dimensional sphere $\Lambda_{\rho}$ in the conformal boundary $\Ein_{n}$
of $\AdS_{n_{1}}$, the so-called limit set. Except in \S~\ref{sub.class} we assume that  $\Lambda_{\rho}$ is acausal.
Our main purpose is to prove the reverse of Theorem 1.1 in \cite{part1}, namely:

\begin{teo}
\label{teo:main}
Let $\Gamma$ be a lattice in $\op{SO}_0(1,n)$. Then, any quasi-Fuchsian representation
in $\op{SO}_0(2,n)$ with acausal limit set is $(\op{SO}_0(2,n), \cY)$-Anosov.
\end{teo}

As in \cite{part1} the proof deeply involves anti-de Sitter geometry, and is based on
the fact that quasi-Fuchsian representations are precisely
holonomy representations of spatially compact, globally hyperbolic 
Lorentzian manifolds locally modelled on $\AdS_{n+1}$.
We will also consider the case where no special hypothesis is made on the group $\Gamma$
(see \S~\ref{sub.gammanolattice}). We give arguments in favor of the idea that Anosov representations 
coincide with ``quasi-Fuchsian'' representations even
if $\Gamma$ is not \textit{a priori} assumed to be isomorphic to a lattice in $\SO_{0}(1,n)$.
A crucial point is the fact that
GHC-regular $\AdS$-spacetimes with acausal limit set admit Cauchy surfaces which are $\op{CAT}(-1)$,
implying that the associated holonomy representation is still Anosov, but now in a weaker, non-differentiable sense.
More generally, in the last {\S} we give a list of groups admitting quasi-Fuchsian
representations into $\SO_{0}(2,n)$ that we expect to be exhaustive.

\section{A criteria for Anosov representations}
\label{sub.def3anosov}
A technical difficulty arising when one wants to prove that a representation is Anosov
is to ensure the exponential decay. The following lemma shows that this feature follows from the compactness
of the ambient manifold $N$ of the Anosov flow and a weaker contraction property somewhere along the orbit.
Compare with \cite[{\S}~5.1.1]{part1}.

\begin{prop}
\label{pro.pasdexp}
Let $\rho: \Gamma \to \SO_0(2,n)$ be a representation, and assume the
existence of continuous maps $\ell^{\pm}_{\rho}: \wt{N} \to \Ein_{n}$ and 
of a continuous and $\Gamma$-equivariant family of Riemannian metrics $g^{p}$ defined in a neighborhood 
of $\ell_\rho^{\pm}(p)$ in $\Ein_{n}$ such that, for all $p$ in $\wt{N}$, there is $t>0$
such that, for all $w$ in $\T_{\ell_\rho^+(p)}\Ein_n$ (respectively in $\T_{\ell_\rho^-(p)}\Ein_n$), we have
${g}^{\tilde{\Phi}^{t}(p)}(w,w) \geq 2{g}^{p}(w,w)$ (respectively 
${g}^{\tilde{\Phi}^{t}(p)}(w,w) \leq {g}^{p}(w,w)/2$). Then, $\rho$
is $(\SO_0(2,n), \cY)$-Anosov.
\end{prop}

\begin{proof}
We simply prove that the hypothesis above imply the exponential decay and exponential expansion
expressed in Remark~5.1.2 in \cite{part1}. Let $\pi_\rho: E_\rho \to N$ be
the flat $\Ein_n$-bundle associated to $\rho$, and let $s^\pm: N \to E_\rho$
be the sections induced by $\ell^\pm_\rho$. 
Since the the family $(g^p)_{(p \in \wt{N})}$ is $\Gamma$-equivariant,
it induces for every $p$ in $N$ a metric $g_\pm^p$ on the fiber $\pi^{-1}(p)$ 
near $s^-(p)$ and $s^+(p)$. Denote by $V^\pm(p)$ the vertical tangent bundle at $s^\pm(p)$.
For every $p$ in $N$ and every $t$ define:

\begin{eqnarray*}
\alpha^-(p,t) & = & \sup_{w \in V^-(p)} \frac{g_-^{{\Phi}^{t}(p)}(w,w)}{g_-^p(w,w)}\\
\alpha^+(p,t) & = & \inf_{w \in V^-(p)} \frac{g_+^{{\Phi}^{t}(p)}(w,w)}{g_+^p(w,w)}
\end{eqnarray*}

Obviously, for $s,t >0$: 

\begin{eqnarray*}
\alpha^-(p, t+s) & \leq & \alpha^-(p, s)\alpha^-(\Phi^s(p),t) \\
\alpha^+(p, t+s) & \geq & \alpha^+(p, s)\alpha^+(\Phi^s(p),t) 
\end{eqnarray*}

By hypothesis, and since $N$ is compact, there is a finite covering $(U_i)_{(1 \leq i \leq k)}$ of $N$,
and a sequence $(T_i)_{(1 \leq i \leq k)}$ such that for any $i$ in $\{ 1, \ldots, k \}$
and any $p$ in $U_i$ we have $\alpha^-(p,T_i) \leq 1/2$. Put $T = \sup\{T_i | 1 \leq i \leq k \}$
and $a = \sup\{ \alpha^-(p,t) | t \in [0, T], p \in N\}$. For any $p$ in $N$ there exist sequences
$(t_j)_{(0 \leq j \leq J)}$ and $(i_j)_{(0 \leq j \leq J)}$ such that $t_0 = 0$, 
$t_{J-1} \leq t \leq t_J$, $\Phi^{t_j}(p)$ lies in $U_{i_j}$ and $t_{j+1} = t_j + T_{i_j}$.
Then:
\begin{eqnarray*}
\alpha^-(p,t) & \leq & \alpha^-(p, T_{i_0})\alpha^-(\Phi^{t_1}(p), T_{i_1}) \ldots \alpha^-(\Phi^{t_{J-1}}(p), t-t_{J-1}) \\
 & \leq & (1/2)^{J-1}a \\
 & \leq & a(1/2)^{t/T-1} 
 \end{eqnarray*}

since $t \leq JT$.
It follows that $\alpha^-(p,t)$ decreases exponentially with $t$. Similarly, $\alpha^+(p,t)$ increases
exponentially with $t$. The proposition follows.
\end{proof}

\section{Dynamical properties}
\label{s.dynamique}
In this {\S} we consider a GHC-regular representation $\rho:\Gamma \to \op{SO}_0(2,n)$
with acausal limit set $\Lambda_\rho$. We don't 
assume that the group $\Gamma$ is isomorphic
to a lattice in $\op{SO}_0(1,n)$.
Let $(\gamma_n)_{(n \in \NN)}$ be a sequence in $\Gamma$ escaping to infinity. It will be convenient
to consider the image sequence $(\rho(\gamma_n))_{(n \in \NN)}$ as a sequence in $\SO_0(2,n+1)$ 
through the inclusion $\SO_0(2,n) \subset \SO_0(2,n+1)$ so that our dynamical study applies 
in $\Ein_{n+1}$, and hence in the $\rho(\Gamma)$-invariant conformal copy
of $\AdS_{n+1}$ inside $\Ein_{n+1}$.

In \cite{franceseinstein} (see also~\cite[\S~7]{primer}), C. Frances 
studied the dynamical behavior in $\Ein_{n+1}$
of $(\rho(\gamma_n))_{(n \in \NN)}$. Up to a subsequence, one the following holds
(we will just afterwards remind ingredients of the proof):

\begin{enumerate}

\item \textit{Balanced distortion:} There are two lightlike geodesics $\Delta^+$, $\Delta^-$ in $\Ein_{n+1}$,
called attracting and repelling photons, and two continuous applications 
$\pi_+: \Ein_{n+1} \setminus \Delta^- \to \Delta^+$
and $\pi_-: \Ein_{n+1} \setminus \Delta^+ \to \Delta^-$ such that: 

-- the fibers of $\pi_{+}$ (respectively $\pi_-$)
are past lightcones $C^-(x)$ of points in $\Delta^-$ 
(respectively of points in $\Delta^+$),

-- for every compact subset $K \subset \Ein_{n+1} \setminus \Delta^-$, the sequence
$\rho(\gamma_n)$ uniformly converges on $K$ to $\pi_+$,

-- for every compact subset $K \subset \Ein_{n+1} \setminus \Delta^+$, the sequence
$\rho(\gamma^{-1}_n)$ uniformly converges on $K$ to $\pi_-$.

\item \textit{Unbalanced distortion:} There are two points $x^+$, $x^-$ in $\Ein_{n+1}$, called
attracting and repelling poles, such that:

-- $\sca{x^+}{x^-} \leq 0$, 

-- for every compact subset $K$ of $\Ein_{n+1}$ contained in
$\Omega^-(x^-):=\{ x \in \Ein_{n+1} /  \sca{x}{x^-} < 0 \}$ (resp.
$\Omega^+(x^-):=\{ x \in \Ein_{n+1} /  \sca{x}{x^-} > 0 \}$) the 
sequence $\rho(\gamma_n)$ uniformly converges on $K$ to the constant map $x^+$
(resp. $(x^+)^\ast$),

-- for every compact subset $K$ of $\Ein_{n+1}$ contained in
$\Omega^-(x^+):=\{ x \in \Ein_{n+1} /  \sca{x}{x^+} < 0 \}$ (resp.
$\Omega^+(x^+):=\{ x \in \Ein_{n+1} /  \sca{x}{x^+} > 0 \}$) the 
sequence $\rho(\gamma_n^{-1})$ uniformly converges on $K$ to the constant map $x^-$
(resp. $(x^-)^\ast$).

\end{enumerate}

\begin{remark}
\label{rk.double}
Our presentation differs from Frances formulation. 
Indeed, we consider the dynamic in $\Ein_{n+1}$, which is the
double covering of the Einstein universe as defined in \cite{franceseinstein}
\ie as the projection of $\cC_{n+1}$ in the projective space
$\PP(\RR^{n+3})$, and not the projection in the sphere of rays $\SS(\RR^{n+3})$.
C. Frances had no need to
distinguish future cones from past cones since they have the same projection in 
$\PP(\RR^{n+3})$. 

A nuisance of the option to consider the double covering is the
non-uniqueness of pairs of attracting/repelling poles. Indeed, the
opposite pair $(-x^+, -x^-)$ is also convenient.
Moreover, if $\sca{x^-}{x^+}=0$, we have four choices $(\pm x^+, \pm x^-)$ of
pairs of attracting/repelling poles. 
\end{remark}

\begin{remark}
Every $\rho(\gamma_n)$ belongs to the subgroup $\SO_0(2,n)$ of $\SO_0(2, n+1)$,
\ie preserves the conformal embedding $\AdS_{n+1} \subset \Ein_{n+1}$ and its boundary
$\partial\AdS_{n+1} \approx \Ein_n$. In that situation, all the limit objects $\Delta^\pm$,
$x^\pm$ involved in the various cases in the description of the asymptotic behavior of $(\rho(\gamma_{n_k}))_{(k \in \NN)}$ are
contained in this boundary. In particular, they avoid $\AdS_{n+1}$.
\end{remark}

The dichotomy balanced/unbalanced is based on the Cartan decomposition of $\SO_0(2,n+1)$. More precisely,
consider the quadratic form $\mathfrak{q}_{2,n+1}:=-4a_1b_1 - 4a_2b_2 + x_1^2 + \ldots + x_{n-1}^2$ on
$\RR^{n+3}$: observe that $(\RR^{n+3}, \mathfrak{q}_{2,n+1})$ and $(\RR^{2,n+1}, \mathrm{q}_{2,n+1})$ are isometric
(the isometry is $(a_1, b_1, a_2, b_2, x_1, \ldots, x_{n-1}) \to  
((a_1+b_1)/2, (a_2+b_2)/2, x_1, \ldots, x_{n-1}, (a_1-b_1)/2, (a_2-b_2)/2)$).
Let $\sA$ be the free abelian subgroup of rank $2$ of $\SO_{0}(2,n+1)$ comprising
elements $a(\lambda, \mu)$ acting diagonally on $\RR^{n+3}$ in the coordinates
$(a_1, a_2, b_1, b_2, x_1, \ldots , x_{n-1})$, so that every $x_i$ is unchanged, the coordinates
$a_1$, $a_2$ are multiplied by $\exp(\lambda)$, $\exp(\mu)$, and the coordinates
$b_1$, $b_2$ are multiplied by $\exp(-\lambda)$, $\exp(-\mu)$. It is a real split Cartan subgroup 
of $\SO_0(2,n+1)$, and we consider the Weyl chamber $\sA^+ \subset \sA$ comprising $a(\lambda, \mu)$
such that $0 \leq \mu \leq \lambda$. The Cartan decomposition Theorem ensures that 
every $\rho(\gamma_n)$ can be written in the form $\rho(\gamma_{n})=\sk_{n} a_{n} \ssl_{n}$
such that $a_n=a(\lambda_n,\mu_n)$ belongs to $\sA^+$ and $\sk_n$, $\ssl_n$ belong to 
the stabilizer $\sK$ of the negative definite $2$-plane 
$\{ a_1=b_1, a_2=b_2, x_1=\ldots=x_{n-1}=0 \}$ (it is a maximal compact subgroup).
Observe that elements of $\sK$ are isometries of the Euclidean norm 
$\Vert (u, v, x_{1}, \ldots, x_{n}) \Vert_{0}^{2}:=u^{2} + v^{2} + x_{1}^{2} + ... + x_{n}^{2}$.
Since $(\rho(\gamma_n))_{(n \in \NN)}$ escapes from any compact,
the sequence $(\lambda_n)_{n \in \NN)}$ is not bounded from above.
By compactness of $\sK$ there is a \textit{converging subsequence,\/} 
\ie a subsequence $(\gamma_{n_{k}})_{(k \in \NN)}$ such that $\sk_{n_k}$, $\ssl_{n_k}$
converge to some elements $\sk_\infty$, $\ssl_\infty$ of $\sK$, and such that
$\lim_{k \to +\infty} \lambda_{n_k} = +\infty$, and $\lim_{k \to +\infty} \exp(\mu_{n_k}-\lambda_{n_k})=\nu$
with $0\leq\nu\leq1$.

\subsubsection{Balanced distortion}
The balanced distortion case is the case $\nu>0$. Denote by $P^-$, $P^+$ the codimension two subspaces
$\{ a_1 = a_2 =0 \}$ and $\{ b_1=b_2=0 \}$ respectively. Consider the following linear endomorphisms of $\RR^{n+3}$:

\begin{eqnarray*}
\Pi^+_0(a_1, b_1, a_2, b_2, x_1, \ldots, x_{n-1}) & = & (a_1, 0, {\nu}a_2, 0, \ldots, 0)\\
\Pi^-_0(a_1, b_1, a_2, b_2, x_1, \ldots, x_{n-1}) & = & (0,b_1, 0,{\nu}b_2, 0, \ldots, 0)
\end{eqnarray*}

They induce maps $\pi^+_0: \SS(\RR^{n+3}) \setminus \SS(P^-)$ and 
$\pi^-_0: \SS(\RR^{n+3}) \setminus \SS(P^+)$. Clearly, as a sequence of
transformations of $\SS(\RR^{n+3})$, 
$(a_{n_k})_{(k \in \NN)}$ converges uniformly on compacts of $\SS(\RR^{n+3}) \setminus \SS(P^-)$
to the map induced by $\pi^+_{0}$, and a similar remark applies for the inverse sequence
$(a^{-1}_{n_k})_{(k \in \NN)}$. It follows that the sequence $(\rho(\gamma_{n_k}))_{(nk \in \NN)}$
converges uniformly on compacts of $\SS(\RR^{n+3})\setminus \SS(\ssl_\infty^{-1}P^-)$ towards
$\sk_\infty\circ\pi_0^+\circ\ssl_\infty$ and that $(\rho(\gamma^{-1}_{n_k}))_{(nk \in \NN)}$
converges uniformly on compacts of $\SS(\RR^{n+3})\setminus \SS(\sk_{\infty}P^+)$ towards
$\ssl^{-1}_\infty\circ\pi_0^-\circ\sk_\infty^{-1}$. The description of the dynamic
in $\Ein_{n+1}$ given above follows by observing that the intersections
$P^\pm \cap \cC_{n+1}$ are isotropic 2-planes.

\subsubsection{Unbalanced distortion}
It is the case $\nu=0$. 
Identify the sphere $\SS(\RR^{n+3})$ of rays with 
the $\Vert_{0}$-unit sphere. 
The attracting fixed points of the action of $a_n$ in $\SS(\RR^{n+3})$ are $\pm x^+_0$
where $x^{+}_{0}=(1, 0, \ldots, 0)$, and the repelling fixed points
are $\pm x^-_0$ where $x^{-}_{0}=(0, 0, 1, 0, \ldots, 0)$.
Observe that the $\mathfrak{q}_{2,n+1}$-orthogonal $(x_{0}^{+})^{\orth}$ is the 
hyperplane $\{ b_{1}=0 \}$: it is also
the orthogonal of $x_{0}^{-}$ for the Euclidean norm $\Vert_{0}$. 
Similarly, $(x_{0}^{-})^{\orth}=\{ b_{2}=0 \}$ is the
$\Vert_0$-orthogonal of $x_{0}^{+}$.

For every $\epsilon >0$
let $C^{-}_{0}(\epsilon)$ be the
spherical ball in $\SS(\RR^{n+3})$ of radius $\pi/2-\epsilon$ centered at $x_{0}^{+}$.
It can also be defined as the connected component containing $x_0^+$ of the complement in $\SS(\RR^{n+3})$ of
the $\epsilon$-neighborhood of $(x_{0}^{-})^{\perp}$.

Every vector in $\RR^{n+3}$ splits as a sum $rx_{0}^{+} + y$ with $y$ in $(x_{0}^{-})^{\orth}$. 
Under the action of $a(\lambda_{n}, \mu_{n})$ the component
$rx_{0}^{+}$ is multiplied by $\exp(\lambda_{n})$ whereas the norm of the component $y$ is multiplied
by at most $\exp(\mu_{n})$. It follows easily:

\begin{lemma}
Let $a(\lambda_{n}, \mu_{n})$ be a sequence in $\sA^{+}$ with no balanced distortion. 
For any $\epsilon > 0$ and any $\eta > 0$ there is $N>0$ such that, for every $n>N$,
the restriction of $a(\lambda_{n}, \mu_{n})$ to $C^{-}_{0}(\epsilon)$ is $\eta$-Lipschitz, with
image contained in $C^{+}_{0}(\pi-\eta)$.\fin
\end{lemma}

The description of the dynamic of unbalanced converging subsequences on $\Ein_{n+1}$ given above 
follows easily; the attracting pole $x^+$ is simply the image of $x_0^+$ by $\sk_\infty$, and
the repelling pole is $x^-=\ssl^{-1}_\infty{x}^-_0$. We entered in such a detail that  
the next lemma is now obvious. Consider the hemisphere $\cD^-=\{ x \in \SS(\RR^{n+3}) / \sca{x}{x^-} < 0 \}$.
For every $\epsilon>0$ let $C^{-}(\epsilon)$ be the set of points in $\cD^{-}$ at distance $\geq \epsilon$
from $(x^{-})^{\orth} \cap \SS(\RR^{n+3})$.
Since $\sk_{n}$, $\ssl_{n}$ are isometries for $\Vert_{0}$:

\begin{lemma}
For any $\epsilon > 0$ and any $\eta > 0$ there is $N>0$ such that, for every $k>N$,
the restriction of $\rho(\gamma_{n_{k}})$ to $C^{-}(\epsilon)$ is $\eta$-Lipschitz,
with image contained in the ball centered at $x^{+}$ and of radius $\eta$.\fin
\end{lemma}

The statement we actually need is:

\begin{cor}
\label{cor.dilate}
Assume that $x^{+}$ belongs to $\rho(\gamma_{n_{k}})\cD^{-}$ for $k$ sufficiently big.
Then the differential at of the inverse of $\rho(\gamma_{n_{k}})$, as a transformation
of the unit sphere $\SS(\RR^{2,n+1})$ expands all
the vectors tangent to the sphere at $x^{+}$ by at least a factor $\nu_{k}$, such that
$\nu_{k} \to +\infty$ when $k \to +\infty$.\fin
\end{cor}

Balanced distortion is the typical behavior of converging subsequences $(\rho(\gamma_{n_{k}}))_{(k \in \NN)}$
when $\rho(\Gamma)$ acts properly discontinuously on $\AdS_{n+1}$. But our situation here
is different: by hypothesis, the group $\rho(\Gamma)$ preserves an achronal limit set $\Lambda_\rho$,
which is not pure lightlike since $E(\Lambda_\rho) \neq \emptyset$ (\cite[Lemma~3.6]{part1}).

\begin{prop}
\label{pro.nobalance}
No sequence in $\rho(\Gamma)$ has balanced distortion.
\end{prop}

\begin{proof}
Assume \textit{a contrario} that some sequence $(\rho(\gamma_n))_{(n \in \NN)}$ 
has balanced distortion. Denote by $\Delta^\pm$ the repelling and attracting photons. 
Since $\Lambda_\rho$ is an acausal topological sphere, it intersects $\Delta^+$ at an unique point $x^+$
Since $\Lambda_\rho$ is $\rho(\Gamma)$-invariant,
the image by $\pi^+$ of $\Lambda_\rho \setminus \Delta^-$ is $x^+$.
The fibers of $\pi^+$ are past cones of elements of $\Delta^-$. Hence, 
$\Lambda_\rho \setminus \Delta^-$ is contained in the past cone $C^-(x^+)$
Since $\Delta^- \cap \Lambda_\rho$ 
is a compact embedded segment, $\Lambda_\rho \setminus \Delta^-$ is dense in $\Lambda_\rho$
(this argument is correct when the dimension of
$\Lambda_\rho$ is $\geq 2$. For the case where $\Lambda_\rho$ is a topological circle,
see \cite{mess1, mess2}, or \cite[\S~6.2]{BBZ}).
Hence $\Lambda_\rho$ is contained in $C^-(x^+)$. It is impossible since $\Lambda_{\rho}$ is not pure lightlike.
\end{proof}

\begin{remark}
\label{rk.xdanslambda}
The ambiguity on the definition of pairs of attracting/repel\-ling poles, mentioned in Remark~\ref{rk.double},
can be removed for GHC-regular representations by selecting as poles the ones contained
in $\Lambda_\rho$. Indeed:

-- $\Lambda_\rho$ contains an attracting pole. Indeed, since it is not contained in a cone $C(\pm x^-)$, 
$\Lambda_\rho$ intersects $\Omega^+(x^-)$ or $\Omega^-(x^-)$, and the $\rho(\gamma_n)$-orbit of a point in this intersection
accumulates on $\pm x^+$, that therefore belongs to $\Lambda_\rho$. Similarly,
$\Lambda_\rho$ contains a repelling pole. 

-- $\Lambda_\rho$ contains one and only one attracting pole. Indeed, $x^+$ and $-x^+$ cannot 
both belong to $\Lambda_\rho$ since $\Lambda_\rho$ is not pure lightlike. Similarly,
$\Lambda_\rho$ contains one and only one repelling pole.

Observe that the condition $\sca{x^-}{x^+} \leq 0$ is fulfilled since $\Lambda_\rho$ is achronal.
\end{remark}

\section{Convex hull of GHC-representations}

\subsection{The convex hull}
\label{sub.convexhull}

According to \cite[Lemma~3.9]{part1} the limit set $\Lambda_\rho$ and the regular domain 
$E(\Lambda_\rho)$ are contained in $U \cup \partial{U}$ where $U$ is an affine domain of $\AdS_{n+1}$.
In particular, it is contained in an affine chard $V$ of $\SS(\RR^{2,n})$. We can consider
the convex hull $\op{Conv}(\Lambda_\rho)$ of $\Lambda_\rho$ in this affine chard. This convex hull 
does not depend on the choice of $V$. Moreover, since $E(\Lambda_\rho)$ is convex, it contains
$\op{Conv}(\Lambda_\rho)$ (cf. \cite[Remark~3.11]{part1}). For more details, see for example \cite{barbtz1}.

Alternatively, we also can define $\op{Conv}(\Lambda_\rho)$ as the projection $\SS(C)$ where
$C$ is the set of barycentric combinations $t_1x_1 + \ldots t_kx_k$ where $t_i$ are positive real numbers
such that $t_1 + \ldots t_k=1$ and $x_i$ elements of $\cC_n \subset \RR^{2,n}$ the projections $\SS(x_i)$
of which belong to $\Lambda_\rho$. 

\begin{lemma}
\label{l:convdansE}
The convex hull $\op{Conv}(\Lambda_\rho)$ is compact; its intersection with
$\partial \AdS_{n+1}$ is $\Lambda_\rho$, and the ``finite part'' 
$\op{Conv}(\Lambda_\rho) \cap \AdS_{n+1}=\op{Conv}(\Lambda_\rho) \setminus \Lambda_\rho$
is contained in $E(\Lambda_\rho)$.Fix a $\rho(\Gamma)$-invariant future oriented timelike vector field $V$ on $E(\Lambda_\rho)$.
\end{lemma}

\begin{proof}
The compactness of $\op{Conv}(\Lambda_\rho)$ arises from the compactness of $\Lambda_\rho$.
Let $x=t_1x_1 + \ldots t_kx_k$ be an element of $\RR^{2,n}$ projecting in $\SS(\RR^{2,n})$ on an
element of $\op{Conv}(\Lambda_\rho)$. For every $y$ such that $\SS(y)$ belongs to $\Lambda_\rho$,
according to \cite[Corollary~2.11]{part1}:

\[ \sca{x}{y}=\sum_{i=1}^{k} t_i \sca{x_i}{y} \leq 0 \]

Moreover, if $\sca{x}{y}$ vanishes, then every $\sca{y}{x_i}$ vanishes. But since $\Lambda_\rho$
is acausal, $\sca{y}{x_i}=0$ implies $y=x_i$: according to \cite[Proposition~3.10]{part1}
$\op{Conv}(\Lambda_\rho) \setminus \Lambda_\rho$ is contained in $E(\Lambda_{\rho})$. The lemma follows
since $E(\Lambda_\rho)$ is contained in $\AdS_{n+1}$.
\end{proof}

\begin{lemma}
\label{le.convempty}
If $\op{Conv}(\Lambda_{\rho})$ has empty interior, then $\rho$ is Fuchsian.
\end{lemma}

\begin{proof}
If $\op{Conv}(\Lambda_{\rho})$ has empty interior, it is contained in a projective hyperplane $\SS(v^{\orth})$.
If $\mathrm{q}_{2,n}(v) > 0$ then $\SS(v^{\orth}) \cap \AdS_{n+1}$ is an isometric,
totally geodesic  embedding of $\AdS_n$.
In a well-chosen conformal chard $\AdS_{n+1} \approx \SS^{1} \times \DD^{n}$ 
this $\AdS$-wall is $\{ x_{n}=0 \}$. It is a contradiction since its closure should contain $\Lambda_{\rho}$
whereas $\Lambda_{\rho}$ is a graph over $\partial\DD^{n}$.
Similarly, if $\mathrm{q}_{2,n}(v)=0$ then $\Lambda_{\rho}$ would be pure lightlike.

Hence, up to renormalization, $v$ lies in $\AdS_{n+1}$. If $v' \neq v$ is another element of $\AdS_{n+1}$
then the intersection $\SS(v^{\orth}) \cap \SS((v')^{\orth}) \cap \AdS_{n+1}$, if not empty, is contained
in a totally geodesic hypersurface in $\SS(v^{\orth}) \cap \AdS_{n+1}$: its closure cannot contain the
topological $(n-1)$-sphere $\Lambda_{\rho}$.
Therefore, $v$ is unique: it is a global fixed point for $\rho(\Gamma)$.
\end{proof}

Since we already know that Fuchsian representations are Anosov (\cite[\S~5.2]{part1}), 
\textit{we assume from now that $\op{Conv}(\Lambda_{\rho})$
has non-empty interior.}
The limit set $\Lambda_\rho$ is the projection of an acausal closed subset 
$\wt{\Lambda}_\rho$ in $\wt\Ein_n$ and $E(\Lambda_\rho)$ is the 1-1 projection of
a domain $\wt{E}(\wt{\Lambda}_\rho)$ in $\wt{\AdS}_{n+1} \approx \RR \times \DD^n$.
Recall that there are two maps $f_\rho^-$, $f_\rho^+$ such that 
$\wt{E}(\wt{\Lambda}_\rho)=\{ (\theta, \mathrm{x}) / f_\rho^-(\mathrm{x}) < \theta < f_\rho^+(\mathrm{x}) \}$ 
(cf. \cite[Remark~3.3]{part1}).

\begin{prop}
\label{pro.F-F+}
The complement of $\Lambda_\rho$ in the boundary $\partial\op{Conv}(\Lambda_\rho)$ has
two connected components. Both are closed edgeless achronal subsets of $\AdS_{n+1}$. More precisely,
in the conformal model their lifting in $\wt{\AdS}_{n+1}$ are graphs 
of 1-Lipschitz maps $F_\rho^+$, $F_\rho^-$ from $\DD^n$ 
into $\RR$ such that 

$$f_\rho^- < F_\rho^- < F_\rho^+ < f_\rho^+$$

\end{prop}

For a similar study when $\Lambda_\rho$ is not necessarily a topological sphere
but in the case $n=2$, see \cite[\S~8.10]{barbtz1}. For the proof of this
proposition, we need a few lemmas.

\begin{lemma}
\label{le.timelikecoupe}
Every timelike geodesic of $\AdS_{n+1}$ intersects $\op{Conv}(\Lambda_\rho)$.
\end{lemma}

\begin{proof}
Let $D$ be a timelike geodesic in $\AdS_{n+1}$. It is contained in a totally geodesic
embedding $A$ of $\AdS_2$ in $\AdS_{n+1}$, and the intersection $A \cap \op{Conv}(\Lambda_\rho)$
contains the convex hull in $A$ of $\Cl(A) \cap \Lambda_\rho$. We are thus reduced
to the (easy) case $n=2$. In that case, $A \setminus D$ has two connected components,
and each of them contains a connected component of $\partial{A}$.
The boundary $\partial A$ has two connected components $l_1$, $l_2$, and 
each of these connected 
components is an inextendible timelike curve in $\Ein_1 \subset \Ein_{n+1}$, which therefore
intersects $\Lambda_\rho$ at an unique point $x_i$. Then, the segment $[x_1, x_2]$ 
intersects $D$.
\end{proof}

\begin{lemma}
\label{l.support}
Support hyperplanes in $\SS(\RR^{2,n})$ to $\op{Conv}(\Lambda_\rho)$ at points
inside $\AdS_{n+1}$ are spacelike.
\end{lemma}

\begin{proof}
Let $x$ be a point in $\AdS_{n+1} \cap \op{Conv}(\Lambda_\rho)$, and let $P$ be a support
(projective) hyperplane at $x$ to $\op{Conv}(\Lambda_\rho)$. This support hyperplane
is a projection $\SS(v^\orth)$for some $v$ in $\RR^{2,n}$. If $\mathrm{q}_{2,n}(v) > 0$, then 
$\SS(v^\orth)$ disconnects any affine domain, in particular, the affine domain 
$V$ containing $E(\Lambda_\rho) \cup \Lambda_\rho$, and it follows easily, since $\Lambda_\rho$
is a topological sphere, that the affine hyperplane
$\SS(v^\orth) \cap V$ disconnects $\Lambda_\rho$. It is a contradiction
since this affine hyperplane is a support hyperplane in $V$ and hence cannot
disconnect the convex hull.

If $\mathrm{q}_{2,n}(v)=0$, then the affine hyperplane $V \cap \SS(v^\orth)$ is tangent to
the hyperboloid $\partial U$ at $\SS(v)$ (up to a slight change of affine domain $V$, we can always assume that
$\SS(v)$ belongs to $V$). If it disconnects $\Lambda_\rho$, we obtain a contradiction as above. 
If not, it means that $\SS(v)$ belongs to $\Lambda_\rho$. Write $x$ as a sum $t_1x_1 + \ldots t_kx_k$ where
$x_i$ belongs to $\Lambda_\rho$: $0=\sca{v}{x}=t_1\sca{v}{x_1} + \ldots t_k \sca{v}{x_k}$. Since every $\sca{v}{x_i}$ 
is a nonpositive number, they all vanish, and it implies that $v=x_i$ for every $i$. Hence $x=v$;
it is a contradiction since $x$ is assumed in $\AdS_{n+1}$.
\end{proof}

\begin{proof}[Proof of Proposition~\ref{pro.F-F+}]
Lift $\op{Conv}(\Lambda_\rho)$ in $\uAdS_{n+1} \approx \RR \times \overline{\DD^n}$
as a subdomain $\op{Conv}(\wt{\Lambda}_\rho)$ in $\wt{E}(\wt{\Lambda}_\rho) \cup \wt{\Lambda}_\rho$.
For every $\mathrm{x}$ in $\DD^n$, the line $\RR \times \{ \mathrm{x} \}$ is a timelike geodesic. According
to Lemma~\ref{le.timelikecoupe} it intersects 
$\op{Conv}(\wt{\Lambda}_\rho)$. Moreover, since this intersection is convex, it contains a geodesic segment 
$[F_\rho^-(\mathrm{x}), F_\rho^+(\mathrm{x})] \times \{ \mathrm{x} \}$. If an element $y$ in 
$]F_\rho^-(\mathrm{x}), F_\rho^+(\mathrm{x})[ \times \{ \mathrm{x} \}$ lies on the boundary
of $\op{Conv}(\wt{\Lambda}_\rho)$, then every support hyperplane to the convex hull at the
projection of this point 
must contain the projection of the timelike segment 
$[F_\rho^-(\mathrm{x}), F_\rho^+(\mathrm{x})] \times \{ \mathrm{x} \}$: it contradicts
Lemma~\ref{l.support}. Therefore, the boundary of $\op{Conv}(\wt{\Lambda}_\rho)$ is the union of the graphs
of $F_\rho^{-}$ and $F_\rho^{+}$. It follows quite easily that these graphs are closed, hence,
$F_\rho^{+}$ and $F_\rho^{-}$ are continuous. 

Consider the closed subset $\cE:=\{ F_\rho^{-}=F_\rho^{+}\}$ in $\DD^n$. For every $\mathrm{x}$ in $\DD^n$,
take a small chard in the Klein model around $(F_\rho^{+}(\mathrm{x}), \mathrm{x})$ such that $F_\rho^{-}$ 
and $F_\rho^{+}$ expresses locally as graphs
of maps from an affine hyperplane into $\RR$. Since $\op{Conv}(\Lambda_\rho)$ is convex,
$F_\rho^{+}$ is convex and $F_\rho^{-}$ is concave. It follows that the coincidence locus $\cE$ is also open.
Since $\DD^n$ is connected, if $\cE$ is not empty we get the equality $F_\rho^{-}=F_\rho^{+}$. It is impossible
since the interior of $E(\Lambda_\rho)$ is not empty. Therefore, according to Lemma~\ref{l:convdansE}:

\[ f_\rho^{-} < F_\rho^{-} < F_\rho^{+} < f_\rho^{+} \]

Finally, for every $\mathrm{x}$ in $\DD^n$, let $\SS(v^\orth)$ be a support hyperplane to
$\op{Conv}(\Lambda_\rho)$ at the projection of $(F_\rho^{+}(\mathrm{x}), \mathrm{x})$. According to Lemma~\ref{l.support}, 
$\SS(v^\orth)$ is a totally geodesic embedding of $\HH^n$. In particular, it lifts as
the graph of a 1-Lipschitz map $\varphi^+_v: \DD^n \to \RR$. One of the region
$\{ (\theta,\mathrm{y}) / \theta > \varphi^+_v(\mathrm{y}) \}$, 
$\{ (\theta,\mathrm{y}) / \theta < \varphi^+_v(\mathrm{y}) \}$ is disjoint from
$\op{Conv}(\wt{\Lambda}_\rho)$, and since $F_\rho^{-}(\mathrm{x}) < F_\rho^{+}(\mathrm{x})=\varphi^+_v(\mathrm{x})$, it is
the former. Hence on $\DD^n$ we have $F_\rho^{+} \leq \varphi_p^+$. But since convex domains are
intersections of half-spaces containing them, we get:

\[ F_\rho^{+}=\min_{v} \varphi^+_v \]

Since every $\varphi^+_v$ is 1-Lipschitz, the same is true for $F_\rho^{+}$. Similarly for $F_\rho^{-}$.
\end{proof}

We denote the components of $\partial\op{Conv}(\Lambda_\rho)$ as
$\wt{S}_\rho^+$, $\wt{S}_\rho^-$. 
Denote by $S^+_\rho$, $S^-_\rho$ their projections in $M=\rho(\Gamma)\backslash{E}(\Lambda_\rho)$.

\begin{lemma}
$S^\pm_\rho$ are Cauchy hypersurfaces in $M$.
\end{lemma}

\begin{proof}
Since $\wt{S}^\pm_\rho$ is homeomorphic to $\RR^n$ the quotient 
$S^\pm_\rho=\wt{S}^\pm_\rho$ is a $K(\Gamma,1)$ space, as $\Gamma\backslash\HH^n$.
The cohomology groups $H^n(S^\pm_\rho, \RR)$ and $H^n(\Gamma\backslash\HH^n, \RR)$ are
therefore isomorphic. Since the later is non zero, the former is non zero: the compactness
of $S^\pm_\rho$, and thus the Lemma, follows.
\end{proof}

As an immediate corollary, we get that 
the projection of $\op{Conv}(\Lambda_\rho) \setminus \Lambda_\rho$
is the compact domain of $M$, bounded by the two disjoint Cauchy hypersurfaces $S^\pm_\rho$.
We denote it $C(M)$, and call it the \textit{convex core} of $M$.

\begin{remark}
It can be easily infered from Lemma~\ref{l.support} that $S^{\pm}_{\rho}$ are furthermore
acausal, \ie that $F_\rho^{\pm}$ are contracting. 
\end{remark}

\begin{remark}
\label{rk.caps}
Let $v$ be an element $v$ such that $\mathrm{q}_{2,n}(v) > 0$.
They are totally geodesic embeddings of $\AdS_{n}$. We call \textit{$\AdS$-wall} the intersections of $\AdS_{n+1}$
with the orthogonal $v^{\orth}$. The \textit{half $\AdS$-spaces} defined by $v$ 
are the domains $\op{H}^{+}(v)=\{ x \in \AdS_{n+1} / \sca{v}{x} \geq 0 \}$ and
$\op{H}^{-}(v)=\{ x \in \AdS_{n+1} / \sca{v}{x} \leq 0 \}=\op{H}^{+}(-v)$.
According to Lemma~\ref{le.timelikecoupe}, the intersection between any $\AdS$-wall $\partial\op{H}(v)$,
and $\Lambda_{\rho}$ is a topological $(n-2)$-sphere. Moreover,
in a suitable conformal chard $\op{H}^{+}(v)$
is the domain $\{ (\theta, x_{1}, \ldots x_{n}) \in \SS^{1} \times \DD^{n} / x_{n} > 0 \}$.
It follows that $\sS(v)=\partial\op{H}^{\pm}(v) \cap \op{Conv}(\Lambda_{\rho})$ is a topological
$n$-dimensional disk, in particular compact, and cuts $\op{Conv}(\Lambda_{\rho})$ in two parts
$\sH^{\pm}(v)=\op{H}^{\pm}(v) \cap \op{Conv}(\Lambda_{\rho})$, that we call
\textit{convex caps.} 
\end{remark}

\subsection{Metric on  the convex hull}

In the sequel we need to define a $\rho(\Gamma)$-metric on $\op{Conv}(\Lambda_{\rho})$.
Since the action is cocompact, all these metrics are quasi-isometric one to the other (see \S~\ref{sub.qi}) and the
choice is not important. However, in order to sustain our argumentation, we choose a specific metric.

Let $\Omega$ be a bounded open domain in $\PP(\RR^{n})$, \ie an open domain contained in an affine chard and such that
the closure $\overline{\Omega}$ in this affine chard is compact. The \textit{Hilbert distance} 
between two points $x$, $y$ 
in $\Omega$ is:

\[ d^{H}(x,y):=\log(a,b,x,y) \]

where $a$, $b$ are the two intersections between $\partial\Omega$ and the projective line containing $x$ and $y$,
and where $(a,b,x,y)$ is the cross-ratio. 
It is a distance function, and the associated metric is proper, geodesic and every projective transformation preserving
$\Omega$ preserves the Hilbert distance of points. Moreover, geodesics are intersections between
projective lines and $\Omega$ (see \cite{kelly}).

The interior of the convex hull $\op{Conv}(\Lambda_{\rho})$ is a bounded open domain, hence admits a
well-defined $\rho(\Gamma)$-invariant Hilbert metric. However, in the sequel we will need metrics defined on 
$\op{Conv}(\Lambda_{\rho}) \setminus \Lambda_{\rho}$ and not only on its interior. Hence we have to enlarge
$\op{Conv}(\Lambda_{\rho})$ to another convex domain, still bounded and $\rho(\Gamma)$-invariant, but containing
the boundaries $\wt{S}^{\pm}_{\rho}$.

A suitable solution is to consider, for $\epsilon >  0$ small enough, the domain
$\op{Conv}(\Lambda_{\rho})_{\epsilon}$ in $\AdS_{n+1}$ made of points $x$ such that every 
causal curve in $\AdS_{n+1}$ joining $x$ to an element
of $\op{Conv}(\Lambda_{\rho})$ is of Lorentzian length $\leq \epsilon$. It follows quite easily from the
compactness of $\rho(\Gamma)\backslash\op{Conv}(\Lambda_{\rho})$ that for $\epsilon$ small enough
$\op{Conv}(\Lambda_{\rho})_{\epsilon}$ is contained in $E(\Lambda_{\rho})$. The proof that $\op{Conv}(\Lambda_{\rho})_{\epsilon}$
is still convex is straightforward, we refer to \cite[Proposition~6.31]{BBZ} for a proof formulated in dimension $2+1$, 
but valid in any dimension.
Observe also that $\op{Conv}(\Lambda_{\rho})_{\epsilon}$ is still bounded, and that its interior contains 
$\op{Conv}(\Lambda_{\rho}) \setminus \Lambda_{\rho}$.

In the sequel, we fix once for all $\epsilon$ and denote by $d^{H}$ the restriction to $\op{Conv}(\Lambda_{\rho})$
of the Hilbert metric of $\op{Conv}(\Lambda_{\rho})_{\epsilon}$.

\begin{remark}
\label{rk.ads=hilbert}
If $]a,b[$ is a spacelike geodesic joining two points in $\Ein_{n}$ then for any $x$, $y$
in $]a,b[$ the $AdS$-length of the piece of geodesic between $x$ and $y$ is $\log(a,b,x,y)$ (see e.g. \cite[Theorem 2.2.1.11]{salein},
it is a generalization of the well-known fact that the Hilbert metric on the Klein model of the hyperbolic space is
isometric to the hyperbolic metric).
It follows that in the case where $a$, $b$ lies on $\Lambda_{\rho}$ this length is the Hilbert distance $d^{H}(x,y)$.
\end{remark}

\subsection{Dirichlet domains}
If $\Gamma$ acts freely and properly discontinuously on a proper complete metric space $X$,
there is a well-known way to construct a fundamental domain of its action: the \textit{Dirichlet
domain} (see \cite[pp 243--245]{ratcliffe}). Here the action we consider does not preserve a 
Riemannian metric, but the construction of Dirichlet domain extends easily in our situation:

\begin{defi}
Fix an element $x_0$ of $\op{Conv}(\Lambda_\rho)$. For every $\gamma$ in $\Gamma$, let
$D(\gamma)$ be the domain $\{ x \in E(\Lambda_\rho) / \; \sca{x}{x_0} > \sca{x}{\rho(\gamma) x_0} \}$
(here we consider $E(\Lambda_\rho)$ as a subset of $\AdS_{n+1} \subset \RR^{2,n}$).
The Dirichlet domain $D(\Gamma)$ is the intersection $\cap_{\gamma \in \Gamma} D(\gamma)$.
\end{defi}

\begin{remark}
Since the quotient $M$ is globally hyperbolic, it admits no closed causal curve.
Therefore, $x_0$ and $\rho(\gamma)x_0$ are not causally related: the $\mathrm{q}_{2,n}$-norm of $(\rho(\gamma)x_0 - x_0)$
is positive. The domain $D(\gamma)$ is the interior of the intersection between $E(\Lambda_{\rho})$
and the half $\AdS$-space $\op{H}^{-}(\rho(\gamma){x}_{0}-x_{0})$.
\end{remark}

\begin{lemma}
\label{le:locallyfinite}
The complements $H(\gamma)=E(\Lambda_\rho) \setminus D(\gamma)$ form a locally finite family of subsets of $E(\Lambda_{\rho})$.
\end{lemma}

\begin{proof}
Assume by contradiction that a compact $K$ of $E(\Lambda_\rho)$ intersects infinitely many $H(\gamma_n)$. 
According to Proposition~\ref{pro.nobalance} and Remark~\ref{rk.xdanslambda}
we can assume, up to a subsequence, that the action induced 
in the Klein model by $\rho(\gamma^{-1}_n)$ 
converges uniformly on $K$ to a point $x^-$ in $\Lambda_\rho$. 

On the other hand, there is a sequence of points $(x_n)_{(n \in \NN)}$ in $K$,
converging to some $x$, and such that for every $n$:

\[ \sca{x_n}{x_0} \leq \sca{x_n}{\rho(\gamma_n)x_0}=\sca{\rho(\gamma^{-1}_n)x_n}{x_0} \]

Since the $\rho(\gamma^{-1}_n)x_n$ has
$\mathrm{q}_{2,n}$-norm $-1$, the convergence in the Klein model towards
$x^-$ means that for some sequence $\lambda_n \to 0$ the
$\lambda_n\rho(\gamma_n^{-1})({x}_n)$ converges to a representant $\hat{x}^-$ in $\cC_n$ of $x^-$.

Hence 

\[
\sca{x_n}{x_0} \leq \frac{1}{\lambda_n}\sca{x_0}{\lambda_n\rho(\gamma^{-1}_n)x_n} \]

The left term converges to $\sca{x}{x_0}$, and since $\frac{1}{\lambda_n}$ converges to $+\infty$ and 
$\sca{x_0}{\lambda_n\rho(\gamma^{-1}_n)x_0}$ converges to the negative number $\sca{x_0}{\hat{x}^-}$,
the right term converges to $-\infty$. Contradiction.
\end{proof}

A first corollary of this lemma is that $D(\Gamma)$ is open, and its closure $\Cl(D(\Gamma))$
is the intersection of the closures of the $D(\gamma)$.

\begin{lemma}
\label{le:recouvre}
The $\Gamma$-iterates of $\Cl(D(\Gamma))$ covers $E(\Lambda_\rho)$, \ie:

\[ E(\Lambda_\rho)=\cup_{\gamma \in \Gamma} \rho(\gamma)\Cl(D(\Gamma)) \]

\end{lemma}

\begin{proof}
Let $x$ be in $E(\Lambda_\rho)$; consider the map $\xi: \Gamma \to \RR$ defined
by $\xi(\gamma)=\sca{x}{\rho(\gamma){x}_0}$. 
If there is a sequence $\gamma_n$ such that $\xi(\gamma_n)$ increases, the argument
used in the proof above with the constant sequence $x_n=x$ leads to a contradiction. 
Hence $\xi$ attains its maximum at some $\gamma_0$, \ie
$\sca{\rho(\gamma_0){x}_0}{x} \geq \sca{\rho(\gamma)x_0}{x}$ for every $\gamma$ in $\Gamma$. It follows
that $\rho(\gamma_0)^{-1}x$ belongs to $\Cl(D(\Gamma))$.
\end{proof}

\begin{lemma}
\label{le:disjoint}
The iterates $\rho(\gamma)D(\Gamma)$ are disjoint one from the other.
\end{lemma}

\begin{proof}
If $x$ lies in $\rho(\gamma)D(\Gamma)$, then for every $h$ in 
$\rho(\Gamma) \setminus \rho(\gamma)$ we have:
$$\sca{x}{\rho(\gamma)x_0} > \sca{x}{hx_0}$$
If moreover $x$ lies in $\rho(\gamma')D(\Gamma)$ with $\rho(\gamma') \neq \rho(\gamma)$
then: 
$$\sca{x}{\rho(\gamma')x_0} > \sca{x}{\rho(\gamma)x_0}$$
We obtain a contradiction
with the above in the case $h=\rho(\gamma')$.
\end{proof}

The two lemmas above proves that $\Cl(D(\Gamma))$ is a fundamental domain for the action of
$\rho(\Gamma)$ on $E(\Lambda_\rho)$. From now we restrict to the intersection
$\Cl(D(\Gamma)) \cap \op{Conv}(\Lambda_\rho)$ and denote it $\bar{D}_{\op{conv}}(\Gamma)$.
Since the quotient $C(M)=\rho(\Gamma)\backslash\op{Conv}(\Lambda_\rho)$ is compact:

\begin{prop}
$\bar{D}_{\op{conv}}(\Gamma)$ is a compact fundamental domain for the action of
$\rho(\Gamma)$ on $\op{Conv}(\Lambda_\rho)$.\fin
\end{prop}

This compactness implies that $\bar{D}_{\op{conv}}(\Gamma)$ is the intersection between
the convex hull and a finite sided convex polyhedron. Hence $\bar{D}_{\op{conv}}(\Gamma)$ itself
is also convex.

\subsection{Quasi-isometry between the group and the convex hull }
\label{sub.qi}
A map $f: X \to X'$ between two metric spaces $(X,d)$, $(X', d')$ 
is a \textit{quasi-isometry} if for some $a>0$, $b>0$ we have $(1/a) d(x,y) -b < d'(f(x), f(y)) < a d(x,y) + b$,
and if moreover any point in $X'$ is at distance at most $b$ from the image of $f$.

According to Lemma~\ref{le:locallyfinite} the set $S$ made of elements $\gamma$ of $\Gamma$
such that $\rho(\gamma)\bar{D}_{\op{conv}}(\Gamma) \cap \bar{D}_{\op{conv}}(\Gamma) \neq \emptyset$ is finite,
and according to Lemma~\ref{le:recouvre}, \ref{le:disjoint} $S$ is a generating set of $\Gamma$.
We consider the Cayley graph $(\Gamma_{S}, d_{S})$, \ie the simplicial metric space admitting as vertices
the elements of $\Gamma$, and such that two vertices $\gamma$, $\gamma'$ are connected by an edge of length $1$
if and only if $\gamma'\gamma^{-1}$ lies in $S$. 

Since $\Gamma$ acts cocompactly on $\op{Conv}(\Lambda_{\rho})$, 
the map:,
$$\hat{\jmath}: (\Gamma_{S}, d_{S}) \to (\op{Conv}(\Lambda_{\rho}), d^{H})$$
associating to any vertex $\gamma$ the element $\rho(\gamma)x_0$ of
$\rho(\gamma)\bar{D}_{\op{conv}}(\Gamma)$ is a quasi-isometry. 

A key feature is that the group $\Gamma$ we consider is (Gromov) \textit{hyperbolic;} for definitions and
properties of hyperbolic spaces or groups, we refer to \cite{gromovflot, ghysharpe}.
By definition, the Gromov boundary of a hyperbolic geodesic space $(X,d)$ is the space of complete geodesic rays modulo 
the equivalence relation identifying two rays staying at bounded distance one from the other. 
Any quasi-isometry
between hyperbolic spaces extends as a homeomorphism between their Gromov boundary: the image
by a quasi-isometry of a geodesic ray is \textit{quasi-geodesic}, \ie a map $c: [0, +\infty[ \to X$
such that, for some $a,b >0$:

\[ 1/a | t-s | - b \leq d(c(t), c(s)) \leq a | s-t | +b \]

Moreover, for every $a, b >0$, there is a constant $D$ such that
for  every $(a,b)$-quasi-geodesic ray $c: [0, +\infty[ \to X$ there is a geodesic ray $c: [0, +\infty[$
such that, for every $t$, the distance $c(t)$ to the image of $c_{0}$ is less than $D$, and the distance of $c_{0}(t)$
to the image of $c$ is less than $D$. We say that \textit{$c$ is at bounded distance $\leq D$ from $c_{0}$.}

It follows that the quasi-isometry between $\Gamma_{S}$ and $\HH^{n}$
extends to a homeomorphism between $\partial\Gamma$ and the conformal sphere $\partial\HH^{n}$.

\begin{prop}
\label{pro.ihomeo}
$\hat{\jmath}$ extends as a homeomorphism $\jmath$ between the Gromov boundary $\partial\Gamma \approx \partial\HH^{n}$ and 
the limit set $\Lambda_{\rho}$.
\end{prop}

\begin{proof}
Let $(\gamma_{n})_{(n \in \NN)}$ be the sequence of vertices of $\Gamma_{S}$ visited by a complete geodesic
ray $r_{0}$ in $(\Gamma_{S}, d_{S})$. According to the above, there is a constant $D \geq 0$ such that
the image $\hat{\jmath}(r_{0})$
is at bounded $d^{H}$-distance $\leq D$ from a geodesic ray in $(\op{Conv}(\Lambda_{\rho}), d^{H})$, \ie a projective
segment $[x, y^{+}[$ where $x$ lies in $\op{Conv}(\Lambda_{\rho})$ and $y^{+}$ an element in 
$\partial\op{Conv}(\Lambda_{\rho})_{\epsilon}$. Since this geodesic ray, of infinite $d^{H}$-length,
is contained in $\op{Conv}(\Lambda_{\rho})$ the limit point 
$y^{+}$ actually lies in $\Lambda_{\rho}$.

On the other hand, according to Proposition~\ref{pro.nobalance} every subsequence of 
$(\gamma_{n})_{(n \in \NN)}$ admits a subsequence $(\gamma_{n_{k}})_{(k \in \NN)}$ with mixed or bounded distortion: there is
an attracting pole $x^{+}$ in $\Lambda_{\rho}$ such that $(\rho(\gamma_{n_{k}}))_{(k \in \NN)}$
converges uniformly on compacts of $E(\Lambda_{\rho})$ to the constant map $x^{+}$. In particular,
$x_{k}=\hat{\jmath}(\gamma_{n_{k}})=\rho(\gamma_{n_{k}})(x_{0})$ converges to $x^{+}$.

If $x^{+} \neq y^{+}$ then $]x^{+}, y^{+}[$ is a complete geodesic in $(\op{Conv}(\Lambda_{\rho}), d^{H})$
of infinite length. Hence there is a complete geodesic $c$ in $\Gamma_{S}$ such that $\hat{\jmath}(c)$ is a quasi-geodesic
at bounded distance from $]x^{+}, y^{+}[$. Therefore the geodesic ray $r_{0}$ alternatively approximates
both ends of $c$: it is a contradiction since these ends are distinct whereas a geodesic ray admits ony one accumulation
point in $\partial\Gamma$.

Therefore $x^{+}=y^{+}$. It follows that $x^{+}$ does not depend on the subsequence, and that $y^{+}$
is the extremity of any $d^{H}$-geodesic ray at bounded distance from $\hat{\jmath}(r_{0})$. 
Hence the map $\jmath: [r_{0}] \in \partial\Gamma  \to x^{+} \in \Lambda_{\rho}$
is well-defined.

We now prove the continuity of $\jmath$. Let $V$ be a neighborhood of $x^{+}$ in $\Lambda_{\rho}$.
Let $U$ be a neighborhood of $x^{+}$ in $\Ein_{n+1}$ disjoint from $x_{0}$, such that $U \cap \Lambda_{\rho} \subset V$
and that $U \cap \AdS_{n+1}$ is convex.
Finally, let $\sH^{+}(v)$ be a convex cap contained in $U$ such that $x^{+}$ is in the interior of the topological disk
$\sH^{+}(v) \cap \partial\op{Conv}(\Lambda_{\rho})$. The geodesic segment $[x_{0}, x^{+}[$
crosses $\sS(v)$ at some point $x_{1}$. Let $x_{2}$ be another point of that segment sufficiently close to $x^{+}$ so that the Hilbert distance between
$x_{2}$ and $\sH^{-}(v)$ is bigger that $2D$, where $D$ is the constant such that for every geodesic ray in $\Gamma_{S}$
there is a $d^{H}$-geodesic in $\op{Conv}(\Lambda_{\rho})$ at uniform distance $D$ from $\hat{\jmath}(r)$.

The point $x_{2}$ is at distance $D$ from an element $\rho(\gamma_k)x_0$ of $\hat{\jmath}(r_{0})$. 
Let now $W$ be the neighborhood of $[r_{0}]$ such that every element $[r]$ of $W$ is represented by a geodesic ray $r$
starting from $id$ and containing $\gamma_{k}$. Then, $\hat{\jmath}(r)$ 
is at bounded distance $D$ from the geodesic segment $[x_{0}, \jmath([r])[$.
Hence $[x_{0}, \jmath([r])[$ contains a point $y_{2}$ at distance $\leq D$ from $\rho(\gamma_k)x_0$,
hence at distance $\leq 2D$ from
$x_{2}$. According to our choice of $x_{2}$, this point $y_{2}$ lies on the same side of the wall $\sS(v)$ than
$x_{2}$, \ie in $\sH^{+}(v)$. Hence $[x_{0}, \jmath([r])[$ crosses $\sS(v)$ before reaching 
$\jmath([r])$. Since $U$ is convex, it follows that $\jmath([r])$ lies in $U$, hence, in $V$. The continuity of
$\hat{\jmath}$ is proved.

If $[r]$ and $[r']$ are two distinct elements in $\partial\Gamma$, there is complete geodesic $c: \RR \to \Gamma$ asymptotic to 
$r$ near $-\infty$ and to $r'$ near $+\infty$. The quasi-geodesic $\hat{\jmath}(c)$ is at bounded 
distance from a geodesic $] \jmath([r]), \jmath([r'])[$
in $\op{Conv}(\Lambda_{\rho})$. It follows that $\jmath([r]) \neq \jmath([r'])$.

Finally, for any $x$ in $\Lambda_{\rho}$, the $d^{H}$-geodesic ray $[x_{0}, x[$ is at bounded distance from the image by $\hat{\jmath}$
of a quasi-geodesic ray in $\Gamma_{S}$, hence from the image by $\hat{\jmath}$ of a geodesic ray. It follows that $\jmath$ is onto.
Since $\partial\Gamma$ is compact, the bijective map $\jmath$ is an homeomorphism. The proposition is proved.
\end{proof}

\begin{remark}
\label{rk.versioncontinue}
It was convenient for the proof above to consider $(\Gamma_{S}, d_{S})$.
But this metric space is quasi-isometric in a $\Gamma$ equivariant way to $\HH^{n}$ and also $\T^{1}\HH^{n}$.
Hence, a corollary of Proposition~\ref{pro.ihomeo} is that any quasi-isometry $\hat{\jmath}_{c}: \T^{1}\HH^{n} \to \op{Conv}(\Lambda_{\rho})$
extends as a homeomorphism $\jmath_{c}$ between the Gromov boundary
$\partial\T^{1}\HH^{n}$ and $\Lambda_{\rho}$.
\end{remark}

\subsection{The geodesic flow of the GHC-regular spacetime}
\label{sec.geodflow}
\begin{defi}
The \textit{non-wandering subset,} denoted $\cN(\Lambda_{\rho})$, is the subset of $\cE^{1}\AdS_{n+1}$
comprising elements $(x,v)$ such that the two extremities $\ell^\pm(x,v)$ lie in $\Lambda_{\rho}$.
The \textit{geodesic flow} on $\cN(\Lambda_{\rho})$ is the flow $\tilde{\phi}^{t}_{\cN}$ such that 
$\tilde{\phi}^{t}_{\cN}(x,v) =(x^{t}, v^{t})$
where $x^{t}$ is the point on the geodesic tangent to $(x,v)$ at distance $t$ (along the geodesic) from $x$,
and $v^{t}$ the vector tangent at $x^{t}$ to this geodesic. 

This definition is $\rho(\Gamma)$-equivariant, we denote by $\cN(\rho)$ the quotient of $\cN(\Lambda_{\rho})$
by $\rho(\Gamma)$ and $\phi^{t}_{\cN}$ the flow on $\cN(\rho)$ induce by $\tilde{\phi}_{\cN}^{t}$.
\end{defi}

The projection $\mathrm{p}(\cN(\Lambda_{\rho}))$ is obviously contained in $\op{Conv}(\Lambda_{\rho})$.
Since $\ell^\pm$ are continuous, and since $\Lambda_{\rho}$, $C(M)$ are compact, the quotient
$\cN(\rho)$ is compact. 


\begin{prop}
\label{pro.cNquasi}
There is a $\Gamma$-equivariant homeomorphism $\mathfrak{f}: \T^{1}\HH^{n} \to \cN(\Lambda_{\rho})$ 
mapping orbits of the geodesic flow $\tilde{\phi}^{t}$ on orbits of $\tilde{\phi}^{t}_{\cN}$.
\end{prop}

\begin{proof}
The orbit space of $\tilde{\phi}^t$ is $\partial\HH^n \times \partial\HH^n \setminus \cD$,
whereas the orbit space of $\tilde{\phi}^t_\cN$ is $\Lambda_\rho \times \Lambda_\rho \setminus \cD$
(where $\cD$ denotes the diagonal in both cases). Moreover, the quotient maps 
$p_\phi: \T^1\HH^n \to \partial\HH^n \times \partial\HH^n \setminus \cD$
and $p_\cN: \cN(\Lambda_\rho) \to \Lambda_\rho \times \Lambda_\rho \setminus \cD$
are locally trivial $\RR$-fibrations. By proposition~\ref{pro.ihomeo}, there is
an equivariant homeomorphism $\jmath \times \jmath$ between the orbit spaces;
the question is to lift this homeomorphism in a $\Gamma$-equivariant way 
to a map $\mathfrak{f}$ so that:

\[ p_\cN \circ \mathfrak{f} = (\jmath \times \jmath) \circ p_\phi \]

The way to perform such a lift is quite well-known.
Take a finite collection $(T_i)_{1\leq i \leq l}$ of small transversals to $\tilde{\phi}^t$
in $\T^1\HH^n$ so that for any $p$ in $\T^1\HH^n$ there is a positive real number $t$ in $]-1, +1[$
such that $\phi^t(p)$ lies on $\gamma T_i$ for some $\gamma$ in $\Gamma$. Observe that such a 
family is locally finite: given $x$, there are only finitely many $\gamma$ fulfilling this condition.
Now, since $p_\cN$ is a fibration, and if the $T_i$ are chosen sufficiently small,
for every $i$, the restriction of $(\jmath \times \jmath) \circ p_\phi$
to $T_i$ lifts to a map $\mathfrak{f}_i: T_i \to \cN(\Lambda_\rho)$ such that, on $T_i$:

\[ p_\cN \circ \mathfrak{f}_i = (\jmath \times \jmath) \circ p_\phi \]

For every $p$ in $\T^1\HH^n$, for every triple $\alpha=(i, \gamma, t_i)$ with $-1\leq t_i \leq 1$ such that 
$\tilde{\phi}^t_i(p)$ lies in $\gamma T_i$ define $x_\alpha(p) = \rho(\gamma)\mathfrak{f}_i(\tilde{\phi}^{t_i}(p))$. 
All these points lie on the same $\tilde{\phi}_\cN$-orbit. Now select a partition of unity 
$(f_i)_{1\leq i\leq l})$ on 
$N = \Gamma\backslash\T^1\HH^n$ subordinate to the covering $(U_i)_{1\leq i\leq l})$ where
$U_i = \{ \phi^t(p) / -1 < t < 1, p \in T_i \}$. It associates to every $x_\alpha$ a weight, namely
the value of $f_i$ at the projection in $N$ of $p$. Define $\mathfrak{f}(p)$ as the barycenter
of $x_\alpha$ with respect to these weights. It defines a continuous $\Gamma$-equivariant map
$\mathfrak{f}$ mapping orbits of $\tilde{\phi}^t$ into orbits of $\tilde{\phi}^t_\cN$. 
Now it follows from the hyperbolicity of $\HH^n$ that a diffusion process along the orbits
transform this map to another map, that we still denote $\mathfrak{f}$, which is 
injective along the orbits (see \cite{ghysdiffu, gromovdiffu}). This map obviously satisfies
the condition $p_\cN \circ \mathfrak{f} = (\jmath \times \jmath) \circ p_\phi$ and is $\Gamma$-equivariant.
It follows that it is injective. An homological argument ensures that it is a homeomorphism.
\end{proof}

We can now improve the content of Proposition~\ref{pro.ihomeo}.:

\begin{prop}
\label{pro.lasuitefinale}
For any complete geodesic ray $[x_{0}, x^{+}[$ in $\op{Conv}(\Lambda_{\rho})$
there is a sequence $(\gamma_{n})_{(n \geq 1)}$ in $\Gamma$ and a convex cap $\sH^{+}$ such that:

\begin{enumerate}

\item  the convex caps $\sH_{n}^{+}:=\rho(\gamma_{n})\sH^{+}$ 
shrink uniformly to $x^{+}$,

\item the repelling pole $x^-$ belongs to 
$\mathfrak{d}^-:=\partial\op{Conv}(\Lambda_\rho) \cap \sH^-$,

\item  the attracting pole $x^{+}$ belongs to every 
$\mathfrak{d}^+_n:=\rho(\gamma_n)\mathfrak{d}^+$ where
$\mathfrak{d}^+:=\partial\op{Conv}(\Lambda_\rho) \cap \sH^+$.

\end{enumerate}
\end{prop}

\begin{proof}
For every $x$ in  $[x_{0}, x^{+}[$, let $v(x)$ be the velocity, \ie the unit vector tangent to $[x_{0}, x^{+}[$ 
and oriented towards $x^{+}$. Since $\cN(\rho)$ is compact, the ${\phi}^{t}_{\cN}$-orbit
of the projection of $(x_{0}, v_{0})$ (where $v_0=v(x_0)$) admits an accumulation
point. Let $(x_{\infty}, v_{\infty})$ be a lifting in $\cN(\Lambda_{\rho})$ of this accumulation
point, and let $\sH^{+}$ be a convex cap such that the wall $\sS$ intersects $[x_{0}, x_{+}[$
and such that the final extremity $\ell^{+}(x_{\infty}, v_{\infty})$ lies in the interior of
$\mathfrak{d}^+$. Fix also a positive real number $\epsilon$, and let $W$ be a small neighborhood of $(x_\infty, v_\infty)$ 
in $\cN(\Lambda_\rho)$ made of points of the form $\tilde{\phi}_\cN^t(y,w)$ where:

-- $-\epsilon < t < \epsilon$,

-- $y$ lies in $\sS$,

-- the tangent vector $w$ points in the direction of $\sH^+$, \ie the final extremity of
the $\tilde{\phi}^t_\cN$-orbit of $(y,w)$ lies in the interior of $\mathfrak{d}^+$.

By construction, 
there is a sequence $(\gamma_n)_{(n \in \NN)}$ and a sequence of points $x_n$ in
$[x_0, x^+[$ converging to the final extremity $x^+$ such that $(x_{n}, v_{n})$
(where $v_{n}:=v(x_{n})$) intersects $\rho(\gamma_{n})W$.
By replacing $\gamma_{n}$ by $\gamma_{n}\gamma^{-1}_{1}$ and $(x_{\infty}, v_{\infty})$,
$\sH^{+}$ and $W$ by their images by $\rho(\gamma_{1})$ we can assume that
$\gamma_{1}$ is trivial and that $x_{1}$ belongs to $W$. Hence $x^{+}=\ell^{+}(x_{1}, v_{1})$
lies in the interior of $\mathfrak{d}^{+}$.
Moreover, $\rho(\gamma_n)^{-1}v_{n}$ points in the direction of $\sH^+$: it follows that
$x^+$ belongs to every $\sH^+_n$, and that $x_0$ belongs to $\sH^-_n$.

Up to a subsequence, we can assume that $(\rho(\gamma_n))_{(n \in \NN)}$
is a converging subsequence with unbalanced distortion. 
The $\rho(\gamma_n)x_\infty$ stay at uniformly bounded distance form $x_n$; it follows that they
converge to $x^+$ and that $x^+$ is the attracting pole of $(\rho(\gamma_n))_{(n \in \NN)}$.
On the other hand, every $\sH^-_n$ contains $x_0$: therefore, these convex caps do not shrink
to a point. The repelling pole $x^-$ lies in $\mathfrak{d}^-$. Hence the positive
convex caps $\sH^{+}_{n}$ shrink to the attracting pole $x^+$. The proposition is proved.
\end{proof}

\subsection{End of the proof of Theorem~\ref{teo:main}}
Let $\ell^{\pm}_{\rho}: \T^{1}\HH^{n} \to \Lambda_{\rho}$ be the composition of $\jmath: \partial\Gamma \to \Lambda_{\rho}$ with
$\ell^{\pm}: \T^{1}\HH^{n} \to \partial\Gamma$: they together define a map 
$(\ell^{+}_{\rho}, \ell^{-}_{\rho}): \T^{1}\HH^{n} \to \cY$. In order to achieve the proof of
the main Theorem we just have to construct the metrics $g^{p}$ satisfying 
the hypothesis of Proposition~\ref{pro.pasdexp}.

Fix a $\rho(\Gamma)$-invariant future oriented timelike vector field $V$ on $E(\Lambda_\rho)$.
For every $x$ in $\op{Conv}(\Lambda_{\rho})$
we simplify the notations by denoting simply $h^{x}$ the metric 
$\tilde{g}^{x, V(x)}$ on $\partial{U}({x})\subset \Ein_{n}$ introduced in \cite[\S~5.2.2]{part1}.
We define $g^{x}$ as the metric $h^{\rho(\gamma){x_{0}}}$ where $\gamma$
is an element of $\Gamma$ such that $\rho(\gamma)\bar{D}_{\op{conv}}(\Gamma)$ contains $x$. 
This family of metrics has a drawback: it is not continuous.

A way to construct a continuous family of metrics is the following:
Let $\varsigma: \T^{1}\HH^{n} \to \wt{S}^{+}_{\rho}$ be the composition of the 
homeomorphism $\mathfrak{f}$ of Proposition~\ref{pro.cNquasi} with the 
projection $\pi: \cN(\Lambda_\rho) \to \op{Conv}(\Lambda_\rho)$: it is a $\Gamma$-invariant homeomorphism.
For $p=(x,v)$ in $\T^{1}\HH^{n}$ define the metric $g^{p}_{0}$ as the metric $h^{\varsigma(p)}$ 
on the open neighborhood $\partial{U}({\varsigma(p)})$
of $\ell^{+}_{\rho}(p)$ and $\ell^{-}_{\rho}(p)$. 
These metrics vary continuously with $p$.

Now the key observation is that to check the expanding property for $g^{p}_{0}$ 
is equivalent to check the same property
for $g^{p}$. Indeed:

\begin{lemma}
\label{le.changermetrique}
For every $\delta > 0$, there is a constant $C_{\delta} > 1$ such that for every $x$ and $y$ 
in $\op{Conv}(\Lambda_{\rho})$
such that $d^{H}(x,y) < \delta$, and for every vector $w$ tangent to $\Ein_{n}$ at a point of $\Lambda_{\rho}$
the following inequalities hold:

\[ C_{\delta}^{-1} h^{y}(w,w) \leq h^{x}(w,w) \leq C_{\delta} h^{y}(w,w) \]

\end{lemma}

\begin{proof}[Sketch of proof]
When $y$ is fixed, for example, $y=x_{0}$, the lemma follows from the 
compactness of the $d^{H}$-ball centered at $x_{0}$ 
and the continuity of $x \to h^{x}$. The general case follows by $\rho(\Gamma)$-equivariance.
\end{proof}

Hence, $g^{p}_{0}$ and $g^{p}$ only differ by a factor $C_{\delta}$ where $\delta$ is the
diameter of $\bar{D}_{\op{conv}}(\Gamma)$. Therefore,
the last step in the proof of Theorem~\ref{teo:main} is:

\begin{prop}
\label{pro.numultiplie2}
Let $p=(x,v)$ be an element of $\T^{1}\HH^{n}$.
Then for every $C>0$, there is a time $t>0$ such that for every
tangent vector $w$ to $\Ein_{n}$ at $\ell^{+}_{\rho}(p)$ the inequality
${g}^{\tilde{\phi}^{t}(p)}(w,w) \geq C{g}^{\tilde{p}}(w,w)$ holds.
\end{prop}

\begin{proof}
Let $r_{0}=[x_{0}, x^{+}[$ be the $\pi$-projection
of the image by $\mathfrak{f}$ of the positive $\tilde{\phi}^{t}$-orbit of $p$. Observe that
$x^{+}=\ell^{+}_{\rho}(p)$.  Let $\sH^+$ be the convex cap and
$(\gamma_{n})_{(n \geq 1)}$ be the sequence 
obeying the conclusion of Proposition~\ref{pro.lasuitefinale}. 

According to Lemma~\ref{le.changermetrique}, it is enough to prove that for 
every $C>0$ there is a positive integer $n$ such that 
the $h^{\gamma_{n}x_{0}}$-norm of any $w$ in $\T_{x_+}\Ein_n$
is bounded from below by its $h^{x_0}$-norm multiplied by $C$.
Since the metrics are $\rho(\Gamma)$-equivariant, we have to prove:

$$h^{x_{0}}(d_{x^{+}}\rho(\gamma_n)^{-1}w, d_{x^{+}}\rho(\gamma_n)^{-1}w) \geq Ch^{x_0}(w, w)$$

This inequality only involves the metric $h^{x_0}$. But since $\Lambda_\rho$ is a compact
subset of $\partial{U}({x_0})$, the $h^{x_0}$-norm of vectors tangent to points in $\Lambda_\rho$
is equivalent to their $\Vert_0$-norm - here by $\Vert_0$ we mean the restriction to 
$\Ein_n$ of the spherical metric on $\SS(\RR^{2,n})$ induced by the Euclidean norm.
Hence, to achieve the proof, we just have to check that Corollary~\ref{cor.dilate} applies,
\ie, with the notations introduced in \S~\ref{sec.geodflow}, that the attracting pole $x^+$ belongs
to $\rho(\gamma_n)D^-$. 

The repelling pole $x^-$ belongs to $\mathfrak{d}^-$ (Item (2) of Proposition~\ref{pro.lasuitefinale}).
Hence, the positive convex cap $\sH^+$ is at positive distance from $(x^-)^\orth$ in the
unit sphere $\SS(\RR^{2,n})$, \ie
is contained in $D_\epsilon$ for $\epsilon$ sufficiently small. Hence $\rho(\gamma_n)D^-$ contains 
$\mathfrak{d}_n^+$. Since $x^+$ lies in $\mathfrak{d}_n^+$ (Item (3) of Proposition~\ref{pro.lasuitefinale}),
we obtain as required that $x^+$ belongs
to $\rho(\gamma_n)D^-$.
\end{proof}

\section{Conclusion}

\subsection{Closure of the set of quasi-Fuchsian representations}
In the Riemannian context, the set of quasi-Fuchsian representations
is not closed. But the situation for quasi-Fuchsian representations in $\SO_0(2,n)$ of
lattices in $\SO_0(1,n)$ is different. Whereas quasi-spheres in $\partial\HH^{n+1}$ 
may degenerate, the limit sets of a sequence of quasi-Fuchsian representations $(\rho_k)_{(k \in \NN)}$
in $\SO_0(2,n)$ always converge, up to a subsequence, to a closed achronal topological
sphere $\Lambda$ in $\Ein_n$, since the space of 1-Lipschitz maps $f: \SS^n \to \SS^1$ is compact. 
It is easy to see that if the representations $\rho_k$ converge to some representation $\rho$, 
then $\Lambda$ is preserved by $\rho(\Gamma)$.

\begin{question}
Is $\Lambda$ acausal?
\end{question}

If this question admits a positive answer, the limit representation $\rho$ is Anosov (faithfullness
and discreteness follow from classical arguments).
In other words, Anosov representations would form an entire component of
$\op{Rep}(\Gamma, \SO_0(2,n))$.

An element in favor of a positive answer is the $(2+1)$-dimensional case:
up to finite index, $\SO_{0}(2,2)$ is isomorphic to $\SO_{0}(1,2) \times \SO_{0}(1,2)$,
and quasi-Fuchsian representations (\ie GHC-regular representations) decomposes
as a pair $(\rho_{L}, \rho_{R})$ of cocompact Fuchsian representations the surface group $\Gamma$
into $\SO_{0}(1,2)$. Since Fuchsian representations form 
a connected component of $\op{Rep}(\Gamma, \SO_{0}(1,2))$,
our assertion follows. Moreover, Einstein space $\Ein_{2}$ is homeorphic to
a double covering of $\PP(\RR^{2}) \times \PP(\RR^{2})$ so that the limit set
is a lifting of the graph of a topological conjugacy between the projective actions of
$\Gamma$ on the projective line induced by $\rho_{L}$ and $\rho_{R}$. 
This topological conjugacy is a homeomorphism, meaning that $\Lambda$ is acausal.
For more details, see~\cite{mess1, barbtz1, barbtz2}.

\subsection{Convex cocompact lattices}
Theorem~\ref{teo:main} extends, \textit{mutatis mutandis,\/} 
to the case where $\Gamma$ is a non elementary convex cocompact subgroup
of $\SO_0(1,n)$, \ie a discrete subgroup acting cocompactly
on the convex hull in $\HH^n$ of its limit set in $\partial\HH^n$ (the non elementary
hypothesis meaning that we require that the cardinal of this limit set is infinite). The definition
of Anosov representation extends in this context by taking as dynamical system $(N, \phi^{t})$ not the entire
$\Gamma\backslash\T^1\HH^n$, but the non-wandering subset
of the geodesic flow in $\Gamma\backslash\T^1\HH^n$: it is not anymore a manifold, but a compact
lamination with a flow (the restriction of the geodesic flow). The
set of $(\SO_0(2,n), \cY)$-Anosov representations is open, and it is still true 
that they correspond to faithfull, discrete representations admitting as limit set a 
closed \textit{acausal} subset in $\Ein_n$, but which now is not a topological sphere.

The main difference is that the associated domains $E(\Lambda_\rho)$ in $\AdS_{n+1}$ are not globally
hyperbolic. However, the action of $\rho(\Gamma)$ on $E(\Lambda_\rho)$ is still free, properly discontinuous
and \textit{strongly causal,\/} \ie the quotient spacetime $\rho(\Gamma)\backslash{E}(\Lambda_\rho)$ is
strongly causal. 
In dimension $2+1$ (when $n=2$) these spacetimes are the so-called \textit{BTZ multi-black holes\/}
(see \cite{BTZ, barbtz2}).

\subsection{Other MGHC spacetimes}
In this paper, we focused on the case where $\Gamma$
is a lattice in $\SO_0(1,n)$. But observe that Theorem~4.7 in \cite{part1} (GHC-spacetimes are GH-regular),
Proposition~\ref{pro.nobalance} (no balanced distortion) and \S~\ref{sub.convexhull} (definition of the convex hull
and the boundary surfaces $S^\pm_\rho$) remains true without this hypothesis.

\subsubsection{GHC-representations with acausal limit set are weakly Anosov}
\label{sub.gammanolattice}
In this {\S} we consider a GHC-regular representation $\rho: \Gamma \to \SO_0(2,n)$, but with
no other assumption on the group $\Gamma$. However we assume that $\Lambda_\rho$ is acausal
so that Lemma~\ref{l:convdansE} holds.

Define the length of Lipschitz curves $c: I \to \wt{S}^\pm_\rho$ as the integral
over $I$ of the Lorentzian norm of the tangent vector (defined everywhere), and
then the distance $\tilde{d}^\pm(x,y)$ between two points $x$, $y$ in $\wt{S}^\pm_\rho$ as the infimum of the length of
Lipschitz curves joining $x$ to $y$. It is not hard to 
see that $\tilde{d}^\pm$ is indeed a distance,
providing to $\wt{S}_\rho^\pm$ a length space structure. 

Observe that $(\wt{S}_\rho^\pm, \tilde{d}^\pm)$
is not in general a Riemannian space, neither Finslerian. However,
this metric structure induces the 
manifold topology on $\wt{S}^\pm_\rho$, which admits a compact
quotient: it is therefore a complete, proper metric space.
By generalized Hopf-Rinow Theorem (\cite[Proposition~I.3.7]{bridson}) $(\wt{S}^\pm, \tilde{d}^\pm_\rho)$
is geodesic: between two points $x$, $y$, there is always a curve joining the two points
realizing the distance.

\begin{prop}
\label{pro.hadamard}
$(\wt{S}^\pm_\rho, \tilde{d}^\pm)$ are complete $\op{CAT}(-1)$ spaces.
\end{prop}

For definition of $\op{CAT}(-1)$ spaces, we refer to \cite[\S~2.1]{bridson}
or \cite{ballmann}.

\begin{proof}
We only consider the upper convex boundary $\wt{S}^-_\rho$, the case of $\wt{S}^+_\rho$ 
is similar (or obtained by reversing the time orientation).
According to the Cartan-Hadamard Theorem (see e.g. \cite[Theorem~4.1]{bridson})
to be a $\op{CAT}(-1)$ space is a local property: since $\wt{S}^-_\rho$ is simply connected,
we just have to prove that every point $x$ admits a neighborhood where the $\tilde{d}^-$
is metric of curvature $\leq -1$ (in the sense of \cite[Definition~II.1.2]{bridson}).

In the Klein model $\wt{S}^-_\rho$ is locally the graph of a convex function 
from an open domain of
$\RR^n$ into $\RR$. More precisely, there is a coordinate system $(t, \bar{x}_1, \ldots, \bar{x}_{n})$,
$-\epsilon < x_i < \epsilon$, $-\eta < t < \eta$ 
on a neighborhood $U$ of $x$ so that:

-- $x$ has coordinates $(0, \ldots, 0)$,

-- $U \cap \wt{S}^-_\rho$ is the graph of a convex map 
$\psi: ]-\epsilon, \epsilon[^n \to ]-\eta, \eta[$,

-- $\{ t=0 \}$ is a  support hyperplane for $\psi$,

-- every tangent vector with negative norm for $-d\eta^2 + dx_1^2 + ... + dx_{n}^2$ has
negative norm for the $\AdS$ metric.

Shrinking $\epsilon$ if necessary, we moreover can assume that the gradient of
$\psi$ has almost everywhere $dx_1^2 + ... + dx_{n}^2$-norm less than 1.
By convolution, we obtain smooth convex maps $\psi_\nu$ which uniformly
converge to $\psi$ when the parameter $\nu >0$ converges to $0$. Moreover,
the norm of their gradient is bounded from above by $1$, it follows that the
graphs $S_\nu$ of $\psi_\nu$ are spacelike. Finally, this uniform convergence implies that
for any Lipschitz curve $c: I \to ]-\epsilon, \epsilon[^n$, the AdS length of 
$s \to (c(s), \psi_\nu(c(s))$ uniformly converges to the $\AdS$-length of
$s \to (c(s), \psi(c(s))$. Hence the graphs $S_\nu$, equipped with their
induced (Riemannian) length metric, converge in the Gromov-Hausdorff topology
to the restriction of $\tilde{d}^-$ to $U \cap \wt{S}_\rho^-$ (cf. \cite[Definition~I.5.33]{bridson}).

We can compute the sectional curvatures of the smooth hypersurfaces $S_\nu$.
Since $\psi_\nu$ is convex, its second fundamental form is positive definite,
and since the ambient $\AdS$ metric has sectionnal curvatures $-1$ the Gauss equation implies 
that $S_\nu$ have sectional curvatures $-1$. They are therefore of curvature $\leq -1$.
The proposition follows since Gromov-Hausdorff limits of length spaces
of curvature $\leq -1$ have curvature $\leq -1$ (\cite[Theorem II.3.9]{bridson}).
\end{proof}

$\op{CAT}(-1)$ spaces enjoy many nice properties. For example, they are hyperbolic 
in the Gromov sense; hence the group $\Gamma$ is Gromov hyperbolic. Furthermore:

\begin{cor}[Proposition II.2.2 in \cite{bridson}]
$(\wt{S}^\pm_\rho, \tilde{d}^\pm_\rho)$ are uniquely geodesic: given two 
points $x$, $y$ there is an unique geodesic joining them.\fin
\end{cor}

We therefore can define the \textit{geodesic flow} of $S^\pm_\rho$, even if
$S^\pm_\rho$ has no unit tangent bundle.

\begin{defi}
Let $\wt{\sG}^\pm_\rho$ denote the space of complete unit speed geodesics of $\wt{S}^\pm_\rho$,
\ie isometries $c: \RR \to \wt{S}^\pm_\rho$, endowed with the topology of uniform
convergence on compact subsets. The geodesic flow $\tilde{\phi}^t_\pm$ is the flow defined
by:

\[ \tilde{\phi}^t_\pm(c)(s)=c(s+t) \]

The group $\rho(\Gamma)$ acts naturally, freely and properly discontinuously on 
$\wt{\sG}^\pm_\rho$. We denote by ${\sG}^\pm_\rho$ the quotient space, and by
$\phi^t_\pm$ the flow on $\sG^\pm_\rho$ induced by $\tilde{\phi}^t_\pm$.
\end{defi}

This flow is not differentiable but \textit{weakly} (or {topologically}) Anosov:
there are two continuous $\Gamma$-invariant foliations $\wt{\cF}_{\pm}^{s}$, $\wt{\cF}_{\pm}^{u}$ 
on $\wt{\sG}^{\pm}_{\rho}$, invariant by the geodesic flow such that for every pair
$p$, $q$ of points in the same leaf of $\wt{\cF}^{s}$ (respectively $\wt{\cF}^{u}$)
there is a real number $t_{0}$ such that the distance between $\tilde{\phi}_{\pm}^{t+t_{0}}(p)$
and $\tilde{\phi}^{t}_{\pm}(q)$ decreases (respectively increases) exponentially with $t$.
This claim follows quite easily from the $\op{CAT}(-1)$ property - it is actually a
general property of Gromov hyperbolic spaces admitting compact quotients: see 
\cite[\S~8.3]{gromovflot}, and for more details, \cite{champetier}, \cite{matheus}.
The fact that the spaces we consider are $\op{CAT}(-1)$ greatly simplifies the definition
of the geodesic flow.

It should be clear to the reader that the methods used in the present paper prove that the
GHC-regular representation $\rho$ satisfies the $(\SO_{0}(2,n))$-Anosov property
as defined in \cite[Remark~5.4]{part1} or appearing as hypothesis in Proposition~\ref{pro.pasdexp} - observe that 
in these formulations the
differential of the flow is not involved. The arguments in \cite[\S~5.3]{part1} still apply for this
non-differentiable version of $(G,Y)$-Anosov property. In other words, we can state that
\textit{GHC-regular representations with acausal limit sets are precisely weakly $(G,Y)$-Anosov representations.}
Moreover, we guess that weakly Anosov representations form an open subset of $\op{Rep}(\Gamma, G)$:
the differentiable setting should be avoided through arguments in \cite{sullivan}. Anyway, for the pair $(\SO_{0}(2,n), \cY)$,
it comes through the discussion above - a representation is GHC-regular with acausal limit set if and only if it is $(\SO_{0}(2,n))$-Anosov -
and the fact that GHC-regular representations form an open domain: using the arguments in \cite{mess1}
one can show that small deformations of holonomy representations of MGHC $\AdS$-spacetimes are still holonomy 
representations of MGHC spacetimes (see also the introduction of \cite{bonsantethese} for more details).

Anyway, we don't discuss or justify furthermore this notion of weaky Anosov representations because
we believe that weakly $(G,Y)$-Anosov representations are (differentially) $(G,Y)$-Anosov in the sense of Labourie.
This statement would be a corollary of a positive answer to the following question:

\begin{question}
Let $\rho: \Gamma \to \SO_{0}(2,n)$ be a GH-regular representation with acausal limit set.
Is there a $\rho(\Gamma)$-invariant \textit{smooth} (\ie $C^{r}$ with $r\geq 3$) convex Cauchy hypersurface?
\end{question}

Indeed, we could replace the Cauchy hypersurfaces $\wt{S}^{\pm}$ in the discussion above
by this smooth convex one, \ie with curvature $\leq -1$, hence, with differentiable Anosov geodesic flow.
Concerning this question, observe that the main task in \cite{BBZ} was to give a positive answer to this question
in dimension $2+1$.

Finally, as before, we can address the question:

\begin{question}
Is the space of (weakly) $(\SO_{0}(2,n), \cY)$-Anosov representations closed?
\end{question}

which, as in the case where $\Gamma$ is a lattice of $\SO_{0}(1,n)$, essentially reduces to the proof
that the limit set of a sequence of $(\SO_{0}(2,n), \cY)$-Anosov representations is acausal.

\subsubsection{Classification of MGHC spacetimes of constant curvature $-1$}
\label{sub.class}

\begin{question}
Let $\rho: \Gamma \to \SO_{0}(2,n)$ a GHC-representation with acausal limit set. Is 
$\Gamma$ isomorphic to a lattice of $\SO_{0}(1,n)$?
\end{question}

A natural way to find a positive answer to this question is to exhibit in the associated
MGHC spacetime a Cauchy hypersurface with constant Gauss curvature $-1$: one of 
the main results of \cite{BBZK} is precisely that such a Cauchy hypersurface exists in
the $(2+1)$-dimensional case. Of course, in this low dimension, this kind of argument is 
sophisticated, since it is only a matter to prove that the genus of the Cauchy surfaces is $\geq 2$,
which can be obtained with more elementary arguments. However, this last idea does not
extends in higher dimension, whereas most part of the content of \cite{BBZK} applies
in any dimension.

Another way to give a positive answer would be to study the functional on $\op{Anos}_{\cY}(\Gamma, \SO_{0}(2,n))$
associating to a representation the volume of the convex core in the associated spacetime.
Indeed, according to Lemma~\ref{le.convempty}, this functional vanishes only on Fuchsian representations. 

Finally, it is easy to produce GHC-regular representations with non-acausal limit set: let $(p,q)$ be a pair of positive integers
such that $p+q=n$, and let $\Gamma$ be a cocompact lattice of $\SO_{0}(1,p) \times \SO_{0}(1,q)$. 
There is a natural inclusion of $\SO_{0}(1,p) \times \SO_{0}(1,q)$ into $\SO_{0}(2,n)$ arising
from the orthogonal splitting $\RR^{2,n}=\RR^{1,p} \oplus \RR^{1,p}$. The isotropic cone of $\RR^{1,p}$ (respectively
$\RR^{1,q}$) is contained in $\cC_{n}$ and projects in $\Ein_{n}$ on the union of two spacelike spheres $\Lambda^{\pm}_{p} \approx \SS^{p-1}$
(respectively $\Lambda^{\pm}_{q} \approx \SS^{q-1}$). Every point in $\Lambda^{\pm}_{p}$ is joined to every point in $\Lambda^{\pm}_{q}$
by a lightlike geodesic segment in $\Ein_{n}$: let $\Lambda$ be the union of lightlike geodesic segments joining a point of
$\Lambda^{+}_{p}$ to a point in $\Lambda^{+}_{q}$ and avoiding $\Lambda^{-}_{p} \cup \Lambda^{-}_{q}$.
We leave to the reader the proofs of the following facts:

-- $\Lambda$ is a non pure lightlike achronal topological sphere,

-- The convex hull $\op{Conv}(\Lambda)$ of $\Lambda$ in $\AdS_{n+1}$ coincide with the regular domain $E(\Lambda)$.

The group $\Gamma \subset \SO_{0}(1,p) \times \SO_{0}(1,q) \subset \SO_{0}(2,n)$ preserves
$\op{Conv}(\Lambda)=E(\Lambda)$; the quotient spacetime $M(\Gamma)=\Gamma\backslash{E}(\Lambda)$
is MGH. Moreover, it is spatially compact: indeed, the set of orthogonal sums $u+v$ where $u$ (respectively $v$)
is an element of $\RR^{1,p}$ such that $\op{q}_{1,p}(u)=-1/2$ (respectively an element of $\RR^{1,q}$
of $\op{q}_{1,q}$-norm $-1/2$) admits two components in $\AdS_{n+1}$, one lying in 
$E(\Lambda)$. This component is a spacelike hypersurface isometric to $\HH^{p} \times \HH^{q}$ and
$\Gamma$-invariant. Its projection in the quotient $M(\Lambda)$ is a compact spacelike hypersurface,
hence a Cauchy hypersurface.

\begin{remark}
By Margulis superrigidity Theorem (\cite{margulis}), if $p,q\geq2$ every $\Gamma$ into $\SO_{0}(2,n)$
either has finite image, or conjugate in $\SO_{0}(2,n)$ to the inclusion 
$\Gamma \subset \SO_{0}(1,p) \times \SO_{0}(1,q) \subset \SO_{0}(2,n)$. 
It follows that every MGHC spacetime of constant curvature $-1$
with fundamental group isomorphic to a lattice $\Gamma$ in $\SO_{0}(1,p) \times \SO_{0}(1,q)$ is isometric
to a spacetime $M(\Gamma)$ described above.
\end{remark}

\begin{remark}
When $n=2$, the only possibility is $p=q=1$. It is the case of \textit{Torus universe}
(see \cite[\S~7]{BBZ}, \cite[\S~3.3]{carlip}).
\end{remark}

\begin{question}
Let $\rho: \Gamma \to \SO_{0}(2,n)$ be a GHC-regular representation with non acausal limit set.
Is $\Gamma$ isomorphic to a lattice of some product $\SO_{0}(1,p) \times \SO_{0}(1,q)$ ?
\end{question}

Our personal guess is that all the questions reported above admit a positive answer. 

\begin{conj}
Every GHC-regular representation into $\SO_{0}(2,n)$ is either a quasi-Fuchsian representation of a lattice in $\SO_{0}(1,n)$, or
a representation of a lattice in $\SO_{0}(1,p) \times \SO_{0}(1,q)$ with $p+q=n$, $p\geq1$, $q\geq1$.
\end{conj}

\bibliography{artads}
\bibliographystyle{alpha}

\end{document}